\documentclass{amsart}
\usepackage{amsfonts}
\usepackage{amsmath}
\usepackage{amsthm}
\usepackage{array}
\usepackage{graphics}
\usepackage[usenames]{color}
\usepackage[a4paper]{geometry}
\usepackage{amssymb}
\usepackage[all]{xypic}
\usepackage{graphicx}

\newtheorem{theorem}{Theorem}

\newtheorem{conjecture}{Conjecture}
\newtheorem{corollary}[theorem]{Corollary}

\newtheorem{remark}[theorem]{Remark}

\newtheorem{lemma}[theorem]{Lemma}

\newtheorem{proposition}[theorem]{Proposition}

\begin{document}

\title{A closed form expression for the Drinfeld modular polynomial $\Phi_T(X,Y)$}
\author{Alp Bassa and Peter Beelen}
%\address{Ecole Polytechnique F\'ed\'erale de Lausanne, EPFL-SFB-IMB-CSAG, Station 8, 1015, Lausanne, Switzerland} %\email{alp.bassa@epfl.ch} % %\author{Peter
%Beelen} %\address{DTU-Mathematics, Technical University of Denmark, Matematiktorvet, Building 303S, DK-2800, Lyngby, Denmark} %\email{p.beelen@mat.dtu.dk} %
%\subjclass{Primary:14H05, Secondary:11R32} % %\keywords{Function field, modular curve, Galois tower}

\thanks{Peter Beelen was partially supported by DNRF (Denmark) and NSFC (China), grant No.11061130539.}
\keywords{Drinfeld Modular polynomial, Catalan Numbers, Drinfeld Modular Curves.}
\subjclass[2000]{Primary: 11F32, 11F52, 05A10, 11B65; Secondary: 11G09, 14G35, 11F03}

\maketitle

\begin{abstract}
In this paper we give a closed form expression for the Drinfeld modular polynomial $\Phi_T(X,Y) \in \mathbb{F}_q(T)[X,Y]$ for arbitrary $q$ and prove a conjecture of Schweizer. A new identity involving the Catalan numbers plays a central role.
\end{abstract}

\section{Introduction}

Let $\mathbb{F}_q$ denote the finite field with $q$ elements. For a polynomial $P(T) \in \mathbb{F}_q[T]$, the Drinfeld modular polynomial $\Phi_P(X,Y)\in \mathbb{F}_q(T)[X,Y]$ defines a model for the Drinfeld modular curve $X_0(P)$. Although the Drinfeld modular polynomials play an analogously fundamental role as the classical modular polynomials, very little is known about their explicit form. For levels $P$ of low degree it is possible to compute $\Phi_P(X,Y)$ explicitly for relatively small values of $q$ only. However, even in the case $P=T$ no explicit expression for $\Phi_T(X,Y)$ is known for general $q$. The main contribution of this paper is to find such an expression for general $q$ by using a new identity involving Catalan numbers. Schweizer \cite{schweizer} studied $\Phi_T(X,Y)$, found an efficient algorithm to compute it for particular values of $q$ and gave two conjectures concerning its structural properties. One of these conjectures was proven in \cite{conj}, the second follows from a closed form expression for the Drinfeld modular polynomial $\Phi_T(X,Y)$ given below.

If $\phi$ is a Drinfeld module of rank two with $j$-invariant $j$ and $\phi'$ is a $P$-isogenous Drinfeld module with $j$-invariant $j'$, then $\Phi_P(j,j')=0$. Moreover, the function field of $X_0(P)$ is given by $\mathbb{F}_q(T)(j,j')$. Since $\phi'$ is also $P$-isogenous to $\phi$ via the dual isogeny, the Drinfeld modular polynomial is symmetric. In case $P(T)$ is an irreducible polynomial, the extension degrees $[\mathbb{F}_q(T)(j,j'):\mathbb{F}_q(T)(j)]$ and $[\mathbb{F}_q(T)(j,j'):\mathbb{F}_q(T)(j')]$ are both equal to $q^{\deg P}+1$.

If $\deg P=1$, without loss of generality we may assume that $P=T$. The genus of the Drinfeld modular curve $X_0(P)$ is zero in this case. Schweizer \cite{schweizer} found an explicit relation between $j$, $j'$ and a uniformizing parameter of the function field of $X_0(T)$. More precisely, he showed that the function field of $X_0(T)$ equals $\mathbb{F}_q(T)(z)$, where
\begin{equation}\label{eq:jandjprime}
j=\dfrac{(z+T)^{q+1}}{z}  \makebox{ and } j'=\dfrac{(z+T^q)^{q+1}}{z^q}.
\end{equation}

This implies that $\Phi_T(X,Y)$ is the unique monic, (symmetric) polynomial of degree $q+1$ in $X$ (and $Y$), such that
\begin{equation}\label{eq:phiT}
\Phi_T\left(\dfrac{(z+T)^{q+1}}{z},\dfrac{(z+T^q)^{q+1}}{z^q}\right)=0.
\end{equation}

As mentioned above, in \cite{schweizer} an algorithm to compute $\Phi_T(X,Y)$ was given. Also a start was made to describe $\Phi_T(X,Y)$ explicitly. More precisely, writing
$$\Phi_T(X,Y)=\sum_{m=0}^{q+1}P_m(Y)X^m$$
and defining $j_0=-T(T^{q-1}-1)^{q+1}$ it was shown that
\begin{itemize}
\item $P_{q+1}(Y)=1$
\item $P_q(Y)=-(Y^q+T(Y-T^q)^{q-1}-T^{q^2}+T^q-T)$
\item $P_1(Y)=-T^{q^2-2q+1}Y(Y-j_0)^{q-1}+(Y+T^{q^2}-T^q)(T^{q^3-q^2}-1)$
\item $P_0(Y)=(Y-j_0)^{q+1}$
\end{itemize}
Moreover two conjectures were formulated concerning the structure of the polynomial $\Phi_T(X,Y)$. The first conjecture concerned a closed form expression for $\Phi_T(X,Y) \bmod{T-1}$, which was proven to hold in \cite{conj}. The second conjecture was the following:
\begin{conjecture}\label{conj:two}
For all integers $m$ satisfying $2\le m \le q-1$ we have $$\dfrac{P_m(Y)}{P_{q+1-m}(Y)}=\left(\dfrac{T^q(Y-j_0)}{Y-T^q}\right)^{q+1-2m}.$$
\end{conjecture}

In this article we will state and prove a closed form expression for $\Phi_T(X,Y)$ (see Theorem \ref{thm:closedformphi}). Conjecture \ref{conj:two} will be a direct consequence of the closed form expression (see Corollary \ref{cor:conj}). In the closed form the so-called Catalan numbers occur. These numbers are defined as \begin{equation}\label{eq:defcat}
C_n:=\frac{1}{n+1}\binom{2n}{n}=\binom{2n}{n}-\binom{2n}{n-1}
\end{equation}
and come up in a variety of combinatorial problems. In the next section we will show an identity involving Catalan numbers which is only valid in characteristic $p>0$. This identity will play an important role in Section \ref{sec:three}, where a proof of the closed formula of $\Phi_T(X,Y)$ is given.

\section{An identity involving the Catalan numbers}

It will be useful to extend the definition of binomial coefficients $\binom{n}{k}$ to integers $n$ and $k$:
\begin{equation}\label{eq:gendef}
\binom{n}{k}=
\left\{
\begin{array}{cl}
\dfrac{n(n-1)\cdots(n-k+1)}{k!} & \makebox{ if $k \ge 0$} \\
\\
0 & \makebox{ otherwise.}
\end{array}
\right.
\end{equation}
The following lemma is a consequences of identities (5.21) and (5.43) in \cite{conc} and will be used later.
\begin{lemma}
Let $r,m,n$, and $\ell$ be integers and assume that $\ell\ge 0$. Then
\begin{equation}\label{eq:trin}
\binom{r}{m}\binom{m}{n}=\binom{r}{n}\binom{r-n}{m-n}
\end{equation}
and
\begin{equation}\label{eq:sum}
\sum_{i=0}^{\ell}(-1)^i \binom{r-i}{\ell}\binom{\ell}{i}=1.
\end{equation}
\end{lemma}
With this lemma, we can prove the following:

\begin{proposition}\label{prop:char0}
Let $e\ge 2$ be a natural number. Then the following equality of polynomials in $\mathbb{Z}[t]$ holds:
\begin{equation}\label{eq:et}
t^{e+1}-(-1)^e=\sum_{n=0}^{\lfloor\frac{e+1}{2}\rfloor}\left( \binom{e-n}{n-1} + \binom{e-n+1}{n} \right)(t-1)^{e+1-2n}t^{n}.
\end{equation}
\end{proposition}
\begin{proof}
Using Newton's binomium to expand $(t-1)^{e+1-2n}$, we can rewrite equation (\ref{eq:et}):
$$
\begin{array}{rcl}
\displaystyle t^{e+1}-(-1)^e & = &
\displaystyle \sum_{n=0}^{\lfloor\frac{e+1}{2}\rfloor}\sum_{m=0}^{e+1-2n}(-1)^{e+1-m} \left( \binom{e-n}{n-1} + \binom{e-n+1}{n} \right)\binom{e+1-2n}{m}t^{n+m}\\
\\
& = &
\displaystyle\sum_{n\ge 0}\sum_{m \ge 0}(-1)^{e+1+m} \left( \binom{e-n}{n-1} + \binom{e-n+1}{n} \right)\binom{e+1-2n}{m}t^{n+m}\\
\\
 & = & \displaystyle\sum_{n\ge 0}\sum_{i\ge n}(-1)^{e+1+i-n} \left( \binom{e-n}{n-1} + \binom{e-n+1}{n} \right)\binom{e+1-2n}{i-n}t^{i}\\
\\
 & = & \displaystyle\sum_{i\ge 0}(-1)^{e+1+i}\sum_{n=0}^i(-1)^n \left( \binom{e-n}{n-1} + \binom{e-n+1}{n} \right)\binom{e+1-2n}{i-n}t^{i}\\
\\
 & = & \displaystyle\sum_{i=0}^{e+1}(-1)^{e+1+i}\sum_{n=0}^{i}(-1)^n \left( \binom{e-n}{n-1} + \binom{e-n+1}{n} \right)\binom{e+1-2n}{i-n}t^{i}.
\end{array}
$$
Finally we see by moving the terms $t^{e+1}$ and $(-1)^e$, that the proposition is equivalent to the identity
\begin{equation}\label{eq:et2}
0=\sum_{i=1}^{e}(-1)^{e+1+i}\sum_{n=0}^{i}(-1)^n \left( \binom{e-n}{n-1} + \binom{e-n+1}{n} \right)\binom{e+1-2n}{i-n}t^{i}.
\end{equation}
At this point, note that by equation (\ref{eq:trin}) we have:
$$\binom{e-n}{n-1}\binom{e+1-2n}{i-n}=\binom{e-n}{i-1}\binom{i-1}{n-1} \ \makebox{ and } \ \binom{e-n+1}{n}\binom{e+1-2n}{i-n}=\binom{e-n+1}{i}\binom{i}{n}.$$
Further, by equation (\ref{eq:sum}), we have for any $i\ge 1$:
$$\sum_{n=0}^i(-1)^n\binom{e-n}{i-1}\binom{i-1}{n-1}=-1 \ \makebox{ and } \ \sum_{n=0}^i(-1)^n\binom{e-n+1}{i}\binom{i}{n}=1.$$
Combining the above, we see that equation (\ref{eq:et2}) and therefore the proposition, holds.
\end{proof}

The following lemma was shown in \cite{conj}:
\begin{lemma}\label{lem:reduce}
Let $q=p^e$ be a power of a prime $p$ and let $n$ and $k$ be integers satisfying $n\ge 0$ and $n-k<q$,. Then
$$\binom{n}{k}\equiv(-1)^{n-k}\binom{q-1-k}{n-k}\pmod{p}.$$
\end{lemma}

This lemma has the following consequence:

\begin{lemma}\label{lem:catalan}
Let $i$ be an integer satisfying $0 \le i<q-1$ and denote by $C_i$ the $i$-th Catalan number. Then
\begin{equation}\label{eq:cnb}
C_i \equiv (-1)^i \left( \binom{q-1-i}{i}+\binom{q-i}{i+1} \right)  \pmod{p}.
\end{equation}
\end{lemma}
\begin{proof}
Using equation (\ref{eq:defcat}) and Lemma \ref{lem:reduce}, we find that
$$C_i=\binom{2i}{i}-\binom{2i}{i-1}\equiv (-1)^i\binom{q-1-i}{i}-(-1)^{i+1}\binom{q-i}{i+1} \pmod{p}.$$
The result now follows.
\end{proof}

We are now ready to state and prove an identity involving Catalan numbers that will turn out to be useful later.
\begin{theorem}
Let $q$ be a power of a prime $p$. Then
\begin{equation}\label{eq:catalan}
\sum_{i=0}^{\lfloor\frac{q-1}{2}\rfloor} C_i \cdot (t(1-t))^{i} \equiv t^{q-1}+(1-t)^{q-1} \pmod{p}.
\end{equation}
\end{theorem}
\begin{proof}
From equation (\ref{eq:et}) we see that (in characteristic zero)
$$\left(\frac{t}{t-1}\right)^{e+1}+\left(\frac{-1}{t-1}\right)^{e+1}=\sum_{i=0}^{\lfloor\frac{e+1}{2}\rfloor} \left(\binom{e-i}{i-1}+\binom{e-i+1}{i}\right) \left(\frac{t}{(t-1)^2}\right)^{i}.$$
Replacing $t$ by $t/(t-1)$ we find
$$t^{e+1}+\left(1-t\right)^{e+1}=\sum_{i=0}^{\lfloor\frac{e+1}{2}\rfloor} (-1)^i \left(\binom{e-i}{i-1}+\binom{e-i+1}{i}\right) \left(t(1-t)\right)^{i}.$$
Choosing $e=q$ and using Lemma \ref{lem:catalan}, we find
$$t^{q+1}+\left(1-t\right)^{q+1} \equiv 1-\sum_{i=1}^{\lfloor\frac{q+1}{2}\rfloor} C_{i-1} \cdot \left(t(1-t)\right)^{i} \pmod{p}.$$
The theorem follows directly from this.
\end{proof}

\begin{remark}
Note that Lemma \ref{lem:catalan} implies that $C_i \equiv 0 \pmod{p}$ for all $i$ such that $(q-1)/2<i<q-1$. This phenomenon has been observed in for example \cite{alku}.
\end{remark}

\begin{remark}
Note that $t$ and $1-t$ are exactly the roots of the polynomial $s^2-s+t(1-t)$. A reformulation of equation (\ref{eq:catalan}) is therefore:
\begin{equation}\label{eq:catalan2}
\sum_{i=0}^{\lfloor\frac{q-1}{2}\rfloor} C_i \cdot t^{i} \equiv s_1^{q-1}+s_2^{q-1} \pmod{p},
\end{equation}
where $s_1$ and $s_2$ are the two roots of the polynomial $s^2-s+t$.
\end{remark}

\begin{remark}
The generating function of the Catalan numbers $c(x):=\sum_{i=0}^\infty C_i x^i$ can be expressed as
$$c(x)=\dfrac{2}{1+\sqrt{1-4x}}.$$ This implies that:
$$c\left(t(1-t)\right)=\dfrac{1}{1-t}=1+t+t^2+\cdots.$$
Since $c(x)$ has integer coefficients, we can reduce this relation modulo a prime number $p$ and deduce that $c\left(t(1-t)\right) \equiv 1/(1-t) \pmod{p}$.
Equation (\ref{eq:catalan}) is a refined version of this, since it can be rewritten as
$$\sum_{i=0}^{\lfloor \frac{q-1}{2} \rfloor} C_i \cdot \left(t(1-t)\right)^i \equiv 1+t+\cdots+t^{q-2}+2t^{q-1} \pmod{p}.$$
Letting $q$ tend to infinity, but holding $p$ fixed, one recovers that $c\left(t(1-t)\right) \equiv 1/(1-t) \pmod{p}$.
\end{remark}

\section{A closed form expression for $\Phi_T(X,Y)$}\label{sec:three}

Now we come to our main result:

\begin{theorem}\label{thm:closedformphi}
Let $q$ be a power of a prime $p$ and define $j_0=-T(T^{q-1}-1)^{q+1}$. Then
\begin{eqnarray}\label{eq:one}
\Phi_T(X,Y) & = & (X+Y-j_0)^{q+1}-XY^q-X^qY+(XY)^q(T^{1-q}-1)+XY(T^{q-1}-1)^{q^2}\notag\\´
 &  & -T^{1-q}XY\sum_{i=0}^{\lfloor\frac{q-1}{2}\rfloor}C_i \cdot (X Y-T^q (X+Y-j_0))^{q-1-2i}(X Y T^{q^2+1})^{i}.
\end{eqnarray}
\end{theorem}

Note that the term $(XY)^qT^{1-q}$ is canceled by a term in the summation corresponding to $i=0$. In fact, after expanding all terms in equation (\ref{eq:one}), one obtains a polynomial in $X$ and $Y$ with coefficients in $\mathbb{F}_q[T]$. Note that $X^{q+1}$ and $Y^{q+1}$ will have coefficient $1$.

\begin{proof}
We will use the characterisation of $\Phi_T(X,Y)$ given in equation (\ref{eq:phiT}) and show that
substituting
\begin{equation}\label{eq:two}
X=\dfrac{(z+T)^{q+1}}{z} \ \makebox{ and } \ Y=\dfrac{(z+T^q)^{q+1}}{z^q}
\end{equation}
on the left-hand side in equation (\ref{eq:one}), one obtains zero.

Equation (\ref{eq:catalan2}) implies that
\begin{equation}\label{eq:three}
\sum_{i=0}^{\lfloor\frac{q-1}{2}\rfloor}C_i \cdot \left(\dfrac{X Y T^{q^2+1}}{(X Y-T^q (X+Y-j_0))^2}\right)^{i}=s_1^{q-1}+s_2^{q-1},
\end{equation}
with $s_1$ and $s_2$ equal to the two solutions to the equation
\begin{equation}\label{eq:four}
s^2-s+\dfrac{X Y T^{q^2+1}}{(X Y-T^q (X+Y-j_0))^2}=0.
\end{equation}
Using equation (\ref{eq:two}), one can show after some computations that
\begin{equation}\label{eq:five}
X Y-T^q (X+Y-j_0) = \dfrac{ (T^{q^2+1}(z+T)^{q-1}+z^{q+1}(z+T^q)^{q-1})(z+T)(z+T^q)  } {z^{q+1}}
\end{equation}
and
\begin{equation}\label{eq:six}
\dfrac{XYT^{q^2+1}}{(X Y-T^q (X+Y-j_0))^2} = \dfrac{T^{q^2+1}z^{q+1}(z+T)^{q-1}(z+T^q)^{q-1}}{ (T^{q^2+1}(z+T)^{q-1}+z^{q+1}(z+T^q)^{q-1})^2 }.
\end{equation}
Equation (\ref{eq:six}) implies that the solutions $s_1$ and $s_2$ to equation (\ref{eq:four}) are given by
\begin{equation}\label{eq:seven}
s_1=\dfrac{T^{q^2+1}(z+T)^{q-1}}{ T^{q^2+1}(z+T)^{q-1}+z^{q+1}(z+T^q)^{q-1} } \makebox{ and } s_2=\dfrac{z^{q+1}(z+T^q)^{q-1}}{ T^{q^2+1}(z+T)^{q-1}+z^{q+1}(z+T^q)^{q-1} }.
\end{equation}
Using equations (\ref{eq:two}), (\ref{eq:three}), (\ref{eq:five}) and (\ref{eq:seven}), we find that
\begin{eqnarray*}
 &   & T^{1-q}XY\sum_{i=0}^{\lfloor\frac{q-1}{2}\rfloor}C_i \cdot (X Y-T^q (X+Y-j_0))^{q-1-2i}(X Y T^{q^2+1})^{i}\\
 & = & \\
 &   & T^{1-q}XY(X Y-T^q (X+Y-j_0))^{q-1}\sum_{i=0}^{\lfloor\frac{q-1}{2}\rfloor}C_i \cdot \left(\dfrac{X Y T^{q^2+1}}{(X Y-T^q (X+Y-j_0))^2}\right)^{i}\\
 & = & \\
 &   & T^{1-q}\dfrac{(z+T)^{2q}(z+T^q)^{2q}}{z^{q^2+q}}\left((T^{q^2+1}(z+T)^{q-1})^{q-1}+(z^{q+1}(z+T^q)^{q-1})^{q-1}\right)\\
 & = & \\
 &   & \dfrac{T^{q^3-q^2}(z+T)^{q^2+1}(z+T^q)^{2q}}{z^{q^2+q}}+\dfrac{(z+T)^{2q}(z+T^q)^{q^2+1}}{T^{q-1}z^{q+1}}.
\end{eqnarray*}
However, a straightforward computation shows that, still with $X$ and $Y$ as in equation (\ref{eq:two}),
\begin{eqnarray*}
 &   & (X+Y-j_0)^{q+1}-XY^q-X^qY+(XY)^q(T^{1-q}-1)+XY(T^{q-1}-1)^{q^2}\\
 & = & \\
 &   & \dfrac{T^{q^3-q^2}(z+T)^{q^2+1}(z+T^q)^{2q}}{z^{q^2+q}}+\dfrac{(z+T)^{2q}(z+T^q)^{q^2+1}}{T^{q-1}z^{q+1}}.
\end{eqnarray*}
This finishes the proof.
\end{proof}

We will now prove Conjecture \ref{conj:two} as a corollary to Theorem \ref{thm:closedformphi}.
\begin{corollary}\label{cor:conj}
Define for $0 \le m \le q+1$ polynomials $P_m(Y)\in \mathbb{F}_q[T,Y]$ by the identity $$\Phi_T(X,Y)=\sum_{m=0}^{q+1}P_m(Y)X^m.$$ Then for all $m$ such that $2 \le m \le q-1$ we have
$$\dfrac{P_m(Y)}{P_{q+1-m}(Y)}=\left(\dfrac{T^q(Y-j_0)}{Y-T^q}\right)^{q+1-2m}.$$
\end{corollary}
\begin{proof}
Note that the expression $(X+Y-j_0)^{q+1}-XY^q-X^qY+(XY)^q(T^{1-q}-1)+XY(T^{q-1}-1)^{q^2}$ only contributes to $P_m(Y)$ if $m\in\{0,1,q,q+1\}$.
Now let $m$ be an arbitrary integer between $2$ and $q-1$. Since $XY-T^q(X+Y-j_0)=X(Y-T^q)-T^q(Y-j_0)$, we see from Theorem \ref{thm:closedformphi} that the coefficient of $X^m$ in $\Phi_T(X,Y)$ is given by
\begin{eqnarray*}
P_m(Y)& = & T^{1-q}Y\sum_{i=0}^{\lfloor \frac{q-1}{2}\rfloor} C_i\cdot \binom{q-1-2i}{m-1-i}(Y-T^q)^{m-1-i} (-T^q(Y-j_0))^{q-m-i} (T^{q^2+1}Y)^{i}\\
 & = & T^{1-q}Y\left( \dfrac{Y-T^q}{-T^q(Y-j_0)}\right)^{m} \cdot \sum_{i=0}^{\lfloor \frac{q-1}{2}\rfloor} C_i\cdot \binom{q-1-2i}{m-1-i} \dfrac{(-T^q(Y-j_0))^{q-i} (T^{q^2+1}Y)^{i}}{(Y-T^q)^{1+i}}.\\
\end{eqnarray*}
Since $$\binom{q-1-2i}{m-1-i}=\binom{q-1-2i}{(q+1-m)-1-i},$$
we see that
$$\dfrac{P_m(Y)}{P_{q+1-m}(Y)}=\dfrac{ \left( (Y-T^q)/(-T^q(Y-j_0)) \right)^{m} }{ \left( (Y-T^q)/(-T^q(Y-j_0)) \right)^{q+1-m} }.$$
The corollary now follows directly.
\end{proof}

\begin{remark}
Using Lemma \ref{lem:catalan} and Theorem \ref{thm:closedformphi}, it is possible to give a closed form expression for $\Phi_T(X,Y)$ not involving the Catalan numbers. However, the resulting expression is not more compact nor easier to prove, as far as the authors know.
\end{remark}

\begin{remark}
Setting $t=x^{q-1}/(1+x^{q-1})\in \mathbb{F}_q(x)$ in equation (\ref{eq:catalan}) one obtains after some manipulations an algebraic dependency between $u=x^q+x\in \mathbb{F}_q(x)$ and $v=x^{q+1}/(x^q+x)\in \mathbb{F}_q(x)$, namely
$$u^{q} \equiv u^{q-1}v^{q}+v+v\sum_{i=0}^{\lfloor\frac{q-1}{2}\rfloor}C_i u^{q-1-i}v^i \pmod{p}.$$ The functions $u$ and $v$ occur in a natural way when defining a certain asymptotically optimal tower \cite{jnt}. This tower was in \cite{elkiesd} shown to have a Drinfeld modular interpretation. This again shows that equation (\ref{eq:catalan}) plays a role in the theory of Drinfeld modular curves. 
\end{remark}

\bibliographystyle{99}

\begin{thebibliography}{A}

\bibitem{alku} R.~Alter and K.K.~Kubota, Prime and prime power divisibility of Catalan numbers, Journal of Combinatorial Theory Series A 15, pp.~243-–256 (1973).

\bibitem{conj} A.~Bassa and P.~Beelen, A proof of a conjecture by Schweizer on the Drinfeld modular polynomial $\Phi_T(X,Y)$, Journal of Number Theory 131, pp.~1276--1285 (2011).

\bibitem{elkiesd} N.D.~Elkies, Explicit towers of Drinfeld modular curves, Progress in Mathematics 202, pp.~189--198 (2001).

\bibitem{jnt} A.~Garcia and H.~Stichtenoth, On the asymptotic behaviour of some towers of function fields over
finite fields, Journal of Number Theory 61, pp.~248--273 (1996).

\bibitem{conc} R.L.~Graham, D.E.~Knuth and O.~Patashnik, Concrete Mathematics, Addison-Wesley Publishing Company, 1989.

\bibitem{schweizer} A.~Schweizer, On the Drinfeld Modular Polynomial $\Phi_T(X,Y)$, Journal of Number Theory 52, pp.~53--68 (1995).

\end{thebibliography}

\end{document}